\documentclass[11pt]{article}

\usepackage{amssymb,amsmath,amsfonts,amsthm}
\usepackage{latexsym}
\usepackage{graphics}
\usepackage{indentfirst}

\newtheorem{thm}{Theorem}
\newtheorem{lem}{Lemma}
\newtheorem{exam}{Example}

\newtheorem{cor}{Corollary}

\def\qed
 {\ifhmode\unskip\nobreak\hfill$\Box$\medskip\fi
 \ifmmode\eqno{\Box}\fi}

\normalsize

\title{{\bf \huge  The number of graphs of given diameter}}
\author{{\bf Zolt\'an F\"uredi}\\
\small R\'enyi Institute of the Hungarian Academy\\[-0.8ex]
\small Budapest, P.O.Box 127, Hungary, H-1364,\\[-0.8ex]
\small \texttt{furedi@renyi.hu}\\[-0.8ex]
\small and\\[-0.8ex]
\small Department of Mathematics\\[-0.8ex]
\small University of Illinois at Urbana-Champaign\\[-0.8ex]
\small Urbana, IL 61801, USA \\[-0.8ex]
\small \texttt{z-furedi@math.uiuc.edu} ${}^{1}$
\and {\bf Younjin Kim}\\
\small Department of Mathematics\\[-0.8ex]
\small University of Illinois at Urbana-Champaign\\[-0.8ex]
\small Urbana, IL 61801, USA \\[-0.8ex]
\small \texttt{ykim36@math.uiuc.edu} ${}^{2}$}
\date{${}$} 

\begin{document}

\maketitle

\renewcommand{\thefootnote}{\empty}
\footnotetext{
 \emph{2000 Mathematics Subject Classification}:
 05C30, 05C80.\hfill    {\tt  [{\jobname}.tex]}\newline
Printed on \today \newline 
${}^1$
Research supported in part by the Hungarian National Science Foundation
 OTKA,  by the National Science Foundation under grant NFS DMS 09-01276,
  and by the European Research Council Advanced Investigators Grant 267195.\\
${}^2$ 
This work was completed while visiting R\'enyi Institute, Budapest, Hungary}

\vskip -1cm

\begin{abstract}
In this paper it is proved that there are constants $0< c_2< c_1$ such that
the number of (labeled) $n$-vertex graphs
 of diameter $d$  is
 $$ (1+o(1)){\frac{d-2}{2}}n_{(d-1)}3^{n-d+1}2^{{n-d+1}\choose{2}}$$
 whenever $n\to \infty$ and $ 3 \leq d \leq n - c_{1} \log n $, where  $n_{(d-1)}=n(n-1)\dots (n-d+2)$.
A typical graph of diameter $d$ consists of a
 combination of an induced path of length $d$ and a highly connected
 block of size $n-d+3$.
In the case  $ d > n- c_{2} \log n $ the
 typical graph has a completely different snakelike structure.
The number of $n$-vertex graphs of diameter $d$
 is $(1+o(1))\frac{1}{2}n_{(d+1)}3^{n-d-1}d^{n-d-1}$
 whenever $n\to \infty$ and $ d > n- c_{2} \log n $.
 \end{abstract}

\section{Introduction, notations}

Let  $\mathcal{G}  ( n, {\rm{diam}} =d)$  be  the class of graphs of
diameter  $d$ on $n$ labeled vertices.
We usually identify the vertex sets with the set of first $n$ integers,
 $[n]=\{1,2,\cdots, n\}$.
It is well known~\cite{B79} that almost all graphs have diameter 2,
$ |\mathcal{G}  ( n, {\rm{diam}} = 2)| = ( 1- o(1)) 2^{n \choose 2} $.
Tomescu~\cite{T96} proved that $|\mathcal{G}  ( n, {\rm{diam}} =d)|=
   2^{n\choose 2}(6\cdot2^{-d}+o(1))^{n}$ for any fixed $d\geq3$ as $n\to \infty$.
Our aim is to give an exact asymptotic and to extend his result for almost
 all $d$ and $n$.

For a graph $G$ and vertex $v$ we use the notation $N(v)$ (or $N_G(v)$) for the
 neighborhood of $v$.
For positive integers $n$ and $k$ we use
$n_{(k)}$ for the $k$-term product $n(n-1)\dots (n-k+1)$.
$\exp_2[x]$ stands for $2^x$ and $\binom{n}{a,b,\dots,z}$ is the multinomial
 coefficient $n!/(a!b!\dots z!)$.

\section{Two classes of diameter $d$ graphs} \label{intro}

Let $S\cup \{ a, b\}$ be an $s+2$-element set, $|S|=s>1$.
Define $\mathcal{H}(S,a,b)$
 as the class of graphs, $G$, with underlying set $S\cup \{ a,b\}$ such that
 the distance between every pair of vertices is at most 2 except for
 $a$ and $b$, their distance is 3.
We have
\begin{equation}\label{eq:11}
 2^{{s}\choose2} 3^{s}\left(1 - c_3 0.9^{s}\right)   <
  |\mathcal{H}(S,a,b)| < 3^{s}2^{{s}\choose2},
  \end{equation}
where $c_3 > 0$ is an absolute constant, independent of $s$.
Indeed, the neighborhoods of $a$ and $b$ are disjoint,
 there are at most $3^s$ possibilities for $(N(a), N(b))$.
This gives the upper bound.
On the other hand, we can get the lower bound by counting the number
 of graphs on $S\cup \{ a, b\}$ with the property that $N(a)\cap N(b)=\emptyset$
 and $N(x)\cap N(y)=\emptyset$ for some $(x,y)\neq (a,b)$, e.g., see Tomescu~\cite{T94}.

\begin{exam} {\bf A block plus a path.}\quad
{\rm Suppose $3 \leq d < n$.
Let $\mathcal{H}_1(n,d)$ be a class graphs of diameter $d$ with vertex sets
 $V:=[n]$  obtained as follows.
Split $V$ into three disjoint non-empty parts $A$, $S$, $B$ with $|A|=i$,
 $|S|=n-d+1$, $|B|=d-1-i$ ($1{\leq}i{\leq}d-2$).
Put a path $(v_0, v_1, \dots, v_{i-1})$ to $A$, a path
 $(v_{i+2}, \dots, v_{d-1}, v_{d})$ to $B$ and a
 copy of $ \mathcal{H}(S,v_{i-1},v_{i+2})$.
} \end{exam}

As the reversed sequences $A'=\{ v_{d}, v_{d-1}$, $\cdots$, $v_{i+2}\}$,
 $B':=\{v_{i-1}$, $\cdots$, $v_{0} \}$ yield the very same graphs,
 we have that the number of graphs in the above class is
\begin{align}\label{eq:2}
 h_1(n,d)\left(1 - c_3 0.9^{n-d}\right) \leq
 |\mathcal{H}_1(n,d)|\leq {\frac{d-2}{2}}n_{(d-1)}3^{n-d+1}2^{{n-d+1}\choose{2}}
 :=h_1(n,d).
\end{align}

\begin{exam} {\bf Snake-like graphs.}\quad
{\rm Suppose $\frac{2}{3}n < d< n$.
Let $ (V_{0},V_{1}, \cdots, V_{d})$ be a partition of $[n]$ into $1$ and $2$ elements
 parts such that $ |V_{0}|= |V_{1}|= |V_{2}|=|V_{d-2}|=|V_{d-1}| = |V_{d}|=1$ and there are no two consecutive
 $2$-element sets (i.e., $|V_{i}|=2$ implies $|V_{i+1}|=1$).
Let's connect each vertex of $V_{i}$ to at least one vertex of $V_{i-1}$, and add
 edges inside the $V_{i}$'s arbitrarily.
The class of graphs obtained this way is denoted by $\mathcal{H}_2(n,d)$. Every $G \in \mathcal{H}_2(n,d)$ is of diameter $d$, and the only pair of vertices of distance $d$ is $\{ V_0, V_d \}$. Let $N_i$ be the set of vertices of $G$ of distance $i$ from $V_d$. We have $N_d= V_0$. If the sequence $N_0, N_1, \dots, N_d$ also satisfies $|N_0|=|N_1|=|N_2|=1$, $|N_{d-2}|=|N_{d-1}|=|N_d|=1$, and $|N_i| \leq 2$  then $G$ appears twice in $\mathcal{H}_2(n,d)$. Denote the class of these graphs by $\mathcal{H}^{2}_2(n,d)$, and let $\mathcal{H}^{1}_2(n,d) = \mathcal{H}_2(n,d) \backslash \mathcal{H}^{2}_2(n,d)$.
}\end{exam}

Every partition gives $2^{n-d-1}3^{n-d-1}$ graphs, and the number of partitions
 is $$
 {n \choose 2}{{n-2} \choose 2}
  \cdots { {n -2(n-d-2)} \choose 2}\times  [ n-2(n-d-1) ]!\times
  {{d-5-(n-d-1)+1} \choose {n-d-1}}.
$$
So this procedure produces $ n_{(d+1)} (2d-3-n)_{(n-d-1)} 3^{n-d-1}$ graphs and the members of $\mathcal{H}^2_2(n,d)$  are counted twice.
Hence
$$
 2|\mathcal{H}^{2}_2(n,d)|+ |\mathcal{H}^{1}_2(n,d)|= n_{(d+1)} (2d-3-n)_{(n-d-1)} 3^{n-d-1}.
  $$
 We have $|N_{1}|=|N_{2}|=|N_{3}|=1$ and $|N_d|=1$. 
 One can see that $\max \{ |N_{d-1}|, |N_{d-2}| \} > 1$ is possible only if $\max \{ |V_{d-3}|, |V_{d-4}| \} = 2$. Similarly, $|N_i| \geq 3$ implies that $|V_{d-i}| = |V_{d-i+2}| = 2$. The number of such partitions $ (V_0, V_1, \dots, V_d)$ is at most $$\frac{n!}{2^{n-d-1}} \times \left( 2{{d-7-(n-d-2)+1} \choose {n-d-2}} + (n-d-2) {{d-7-(n-d-2)+1} \choose {n-d-2}} \right).$$
   The sum in the parentheses is at most $$ (n-d){{d-6-(n-d-1)+1} \choose {n-d-2}}  = (n-d)\times \frac{(n-d-1)}{d-5-(n-d-1) +1} {{d-5-(n-d-1)+1} \choose {n-d-1}}.$$
We obtain $$ 2 |\mathcal{H}_2(n,d)| \leq n_{(d+1)} (2d-3-n)_{(n-d-1)} 3^{n-d-1} \left( 1 + \frac{(n-d)(n-d-1)}{(2d-n-3)}\right).$$

Since
$$ d^{n-d-1}\left( 1-\frac{2(n-d+1)}{d} \right)^{n-d-1}<   (2d-3-n)_{(n-d-1)}\left( 1 + \frac{(n-d)(n-d-1)}{(2d-n-3)}\right) \leq
 d^{n-d-1}, $$
we get for some $c_4>0$
\begin{equation}\label{eq:3}
 (1-c_4\frac{(n-d-1)^2}{n}) h_2(n,d)
  < |\mathcal{H}_2(n,d)| < \frac{1}{2}n_{(d+1)} d^{n-d-1} 3^{n-d-1}:=h_2(n,d).
  \end{equation}

The estimates (\ref{eq:2}) and (\ref{eq:3}) give the lower bounds for
 the next two Theorems.

\section{Results}

\begin{thm}\label{thm1}
There is a constant $c_{1}>0$ such that the following holds.
If $3\leq d < n-c_{1}\log n$ and $n\to \infty$ then
 almost all $n$-vertex graphs of diameter at least $d$ belong to
 $ \mathcal{H}_1(n,d)$, hence
$$
 |\mathcal{G}  ( n, {\rm{diam}} =d)|=
  (1+o(1)){\frac{d-2}{2}}n_{(d-1)}3^{n-d+1}2^{{n-d+1}\choose{2}} . $$
\end{thm}

\begin{thm}\label{thm2}
There exist a constant $c_{2}>0$ such that
for $n- c_{2} \log n < d<  n$, $n \to \infty$ almost all $n$-vertex graphs
 of diameter at  least $d$ belong to $\mathcal{H}_2(n,d)$, hence
$$
 |\mathcal{G} ( n, {\rm{diam}} =d)|= (1+o(1))\frac{1}{2}n_{(d+1)}
   d^{n-d-1}3^{n-d-1}. $$
\end{thm}

\begin{cor}
For $2 \leq d < n-c_1\log n$ or $n > d > n-c_2 \log n$
\begin{equation} \label{eq:cor}
 \lim_{n\rightarrow\infty} \frac{|\mathcal{G}  ( n, {\rm{diam}} \geq  d+1)|}{|\mathcal{G}  ( n, {\rm{diam}} =d)|}=0
\end{equation}
\end{cor}

The equation (\ref{eq:cor})  was proved by Tomescu~\cite{T96}
 for every fixed $d \geq 2$ and 
 by Grable~\cite{B}
for all $2\leq d \ll \sqrt{n}/\log n $. The main ideas of our proofs are rather straightforward, but one needs very careful estimates and calculations.

\section{Lemmas for the upper bound}\label{lemma} 

Let $V$ be an $n$-element set, $x_0\in V$,  and let
  $P:= (N_0,  N_{1}, \cdots, N_{d})$ be an ordered partition of $V$ into
$d+1$ non-empty parts, $N_0=\{ x_0\}$,  $n_i:=|N_i|$.
Let  $\mathcal{G} ( x_{0}, N_{1}, \cdots, N_{d})$  be  the class of graphs $G$
 with vertex set $V$ such that $N_i$ is the $i$'th neighborhood of $x_0$,
 $N_{i}=\{ y \in V : d_G(x_{0},y)=i\}$.
The number of graphs in each partite set is $2^{n_{i}\choose 2}$ and the
 number of bipartite graphs between $N_{i}$ and $N_{i+1}$ with no isolated
 vertex in $N_{i+1}$ is $(2^{n_{i}}-1)^{n_{i+1}}$.
We obtained
\begin{equation}\label{eq:1}
  \left|\mathcal{G} ( x_{0}, N_{1}, \cdots, N_{d})\right|=
   2^{\sum_{i=1}^{d}{{n_{i}} \choose 2}}\prod_{i=1}^{d-1}(2^{n_{i}}-1)^{n_{i+1}}.
  \end{equation}

Taking all possible $(d+1)$-partitions $(  x_0, N_1, \dots, N_d)$ we count
  each graph from $ \mathcal{G}  ( n, {\rm{diam}} =d) $ at least twice.
We have
\begin{align} \label{eq:23}
2|\mathcal{G}  ( n, {\rm{diam}} =d)|
 & \leq 
 \sum_{n_{1}+n_{2}+\cdots+n_{d}= n -1 \atop n_{1},n_{2},\cdots,n_{d}\geq1 }
    {{n} \choose {1, n_{1},n_{2},\cdots,n_{d}}}  2^{\sum_{i=1}^{d}{{n_{i}} \choose 2}}
    \prod_{i=1}^{d-1}(2^{n_{i}}-1)^{n_{i+1}}.
\end{align}
In the rest of the proof we give sharp upper bounds for the right  hand side
 of (\ref{eq:23}).
We will use the following estimate.
\begin{align}
 {}&  {{n} \choose {1, n_{1},n_{2},\cdots,n_{d}}}
  2^{\sum_{i=1}^{d}{{n_{i}} \choose 2}}
    \prod_{i=1}^{d-1}(2^{n_{i}}-1)^{n_{i+1}} \nonumber\\
 &{}\enskip = n_{(d+1)} {{n-d-1} \choose {n_{1}-1,n_{2}-1,\cdots,n_{d}-1}}
    2^{\sum_{i=1}^{d}{{n_{i}} \choose 2}}
    \prod_{i=1}^{d-1}\frac{1}{n_i}(2^{n_{i}}-1)^{n_{i+1}} \nonumber\\
  &{}\enskip \leq  n_{(d+1)}
 {{n-d-1} \choose {n_{1}-1,n_{2}-1,\cdots,n_{d}-1}} \times
  \exp_2\left[ \sum_{1\leq i\leq d} {n_{i} \choose 2} +
    \sum_{1\leq i\leq d-1} (n_{i}n_{i+1}-1)\right].  \label{eq:6}
\end{align}

Define
$$
 f(x_1, \dots, x_d):= \sum_{1\leq i\leq d} \frac{1}{2}x_i^2 +
   \sum_{1\leq i\leq d-1} x_ix_{i+1}.
  $$

\begin{lem}\label{lm1}
Let $x_1, \dots, x_d\geq 0$ be real numbers,
 $\sum_{i}x_i=s$, $m:=\max_{1< i< d} (x_{i-1}+x_i+x_{i+1})$.
Then
\begin{equation}\label{eq:m_is_big}
   f(\mathbf{x})\leq \frac{1}{2}m^2+ \frac{1}{2}(s-m)^2,
  \end{equation}
and
\begin{equation}\label{eq:m_is_small}
   f(\mathbf{x})\leq \frac{3}{4}ms.
  \end{equation}
\end{lem}

\noindent
Proof:\quad
Suppose that $m=x_{k-1}+x_k+x_{k+1}$,
 then $x_{k-2}\leq x_{k+1}$ and $x_{k-1}\geq x_{k+2}$.
We have
\begin{align*}
  f(\mathbf{x})& \leq \frac{1}{2}\left( (\sum x_i)-(x_{k-1}+x_k+x_{k+1})\right)^2
 + \frac{1}{2}(x_{k-1}+x_k+x_{k+1})^2\\
  {}&\quad\quad\quad\quad\quad \quad\quad
 +x_{k-2}x_{k-1}+x_{k+1}x_{k+2}-x_{k-1}x_{k+1} -x_{k-2}x_{k+2}\\
  & = \frac{1}{2}(s-m)^2+\frac{1}{2}m^2+(x_{k-2}-x_{k+1})(x_{k-1}-x_{k+2}).
\end{align*}
Here the last term is non-positive and we get (\ref{eq:m_is_big}).

To show (\ref{eq:m_is_small}) consider
\begin{align*}
  4f(\mathbf{x})+ \sum x_i^2&=
  x_1^2 + (x_1+x_2)^2 + (x_1+x_2+x_3)^2+ \dots + (x_{i-1}+x_i+x_{i+1})^2+\dots \\
   {}& \quad\quad\quad \quad\quad
   \dots + (x_{d-2}+x_{d-1}+x_d)^2+ (x_{d-1}+x_d)^2+x_d^2
   - 2\sum x_ix_{i+2}\\
 {}&\leq m\left( x_1 + (x_1+x_2)+ \dots + (x_{i-1}+x_i+x_{i+1})+  \dots
 + (x_{d-1}+x_d)+x_d\right)\\
&= 3ms.  \quad\quad\quad \quad\quad \quad\quad\quad \quad\quad \quad
  \quad     \quad    \quad \quad\quad\quad \quad\quad
  \quad\quad\quad
  \quad\quad\quad \quad\quad\square
  \end{align*}

We will use this Lemma to bound in the following form
\begin{equation}\label{eq:gf}
 \sum \binom{n_i}{2} + \sum (n_i n_{i+1}-1)=
    f(x_1, \dots, x_d) +\frac{5s}{2} -x_1 -x_d.
  \end{equation}
where $x_i:=n_i-1$ ($1\leq i\leq d$), $s=\sum x_i= n-d-1$.

\section{Proof of the upper bound for Theorem \ref{thm1}}\label{upper} 

From now on, we suppose that $ 3 \leq d < n-c \log n$, where $c$ is a
 sufficiently large constant.
We put the terms of the right hand side of (\ref{eq:23}) into four groups
  according to the relation of $s:=n-d-1$ and $m:=
  \max_{1< i< d}(n_{i-1}+n_i+n_{i+1}-3)$.
 \newline
${}$\quad -- Case 1: $m< 0.6s$,  \newline
${}$\quad -- Case 2: $0.6s \leq m < s-1$, \newline
${}$\quad -- Case 3: $m=s-1$,  \newline
${}$\quad This means that for some $1< i< d$ one has $n_{i-1}+n_i+n_{i+1}=s+2$,
  there is an $n_t=2$ ($t\neq i-1, i,i+1$) and all other $n_j=1$.
  We consider three subcases
  \newline
${}$\quad -- -- Case 3.1: $t\neq i-2, i+2$,  \newline
${}$\quad -- -- Case 3.2: $t=i-2$, $n_{i+1}\geq 2$, \newline
${}$\quad -- -- Case 3.3: $t=i+2$, $n_{i-1}\geq 3$, \newline
${}$\quad -- Case 4: $m=s$.\newline
${}$\quad We have $n_{i-1}+n_i+n_{i+1}=s+3$, all other $n_j=1$.
Again we handle three subcases separately \newline
${}$\quad -- -- Case 4.1: $n_{i-1}\geq 2$, $n_{i+1}\geq 2$,  \newline
${}$\quad -- -- Case 4.2: $n_0=n_1=\dots n_{d-2}=1$,  $n_{d-1}+n_d=s+2$, \newline
${}$\quad -- -- Case 4.3: $n_{i-1}+n_i=s+2$, $1< i<d$, all other $n_j=1$. \newline
These exhaust all possibilities.
We will show that the sum in each of the above group is $o(h(n,d))$, except
 in the Case 4.3.
We denote by $\Sigma_1$, $\Sigma_2$, $\Sigma_{31}, \dots$ the sum of the right hand side
 of (\ref{eq:1}) corresponding to the above cases.

\medskip

{\bf Case 1.}\quad
To get an upper bound we use (\ref{eq:6}), rearrange, and then (\ref{eq:gf}) and
 finally (\ref{eq:m_is_small}). We have
\begin{align}
 &{}\quad \Sigma_1:= \sum_{{n_{1}+n_{2}+\cdots+n_{d}=n-1 \atop n_{1},n_{2},\cdots,n_{d}\geq 1}
  \atop {m < 0.6s}}
    {{n} \choose {1, n_{1},n_{2},\cdots,n_{d}}}  2^{\sum_{i=1}^{d}{{n_{i}} \choose 2}}
    \prod_{i=1}^{d-1}(2^{n_{i}}-1)^{n_{i+1}} \nonumber\\
 &\leq
 n_{(d+1)}
  \sum_ {m < 0.6s} \left( {{n-d-1} \choose {n_{1}-1,n_{2}-1,\cdots,n_{d}-1}}
 \times
  \exp_2 \left[\sum \binom{n_i}{2} + \sum (n_i n_{i+1}-1)\right] \right)\nonumber\\
 &\leq
 n_{(d+1)}
 \left( \sum {{n-d-1} \choose {n_{1}-1,n_{2}-1,\cdots,n_{d}-1}}
 \right) \nonumber\\
  &{} \quad\quad\quad\quad\quad\quad\quad\quad\quad\quad\quad\quad\times
  \exp_2 \left[\max_{m<0.6s} \left\{ \sum \binom{n_i}{2} + \sum (n_i n_{i+1}-1)
 \right\}\right] \nonumber\\
&\leq  n_{(d+1)}d^{n-d-1}\times
   \exp_2 \left[\max_{m<0.6s}  f(\mathbf{x}) +\frac{5s}{2}\right] \nonumber\\
 &= n_{(d+1)}d^{n-d-1} \exp_2[(3/4) (0.6s)s +5s/2].   \label{eq:10}
\end{align}
This implies
\begin{equation*} 
  \log_2 \Sigma_1 \leq  \log (n_{(d+1)}) +s \log d+0.45s^2 +2.5s.
  \end{equation*}
On the other hand (\ref{eq:2}) gives
\begin{equation}\label{eq:log_hnd}
  \log_2 h_1(n,d)= -1+ \log(d-2)+ \log (n_{(d-1)}) +  (s+2)\log 3 + \binom{s+2}{2}.   
  \end{equation}
A little algebra gives $\log h_1(n,d)-\log \Sigma_1 > s^2/20 -s\log d$
 (for $n-d-1> 100$) and this goes to infinity as $s\to \infty$
 because  $d< n- 41 \log n$ implies $n-d-1> 40\log n> 40 \log d$.
Thus in this range $\Sigma_1 = o(h_1(n,d))$.
\medskip

{\bf Case 2.}\quad
To get an upper bound we use (\ref{eq:6}) but rearrange more carefully.
We have
\begin{align}
 &{}\quad \Sigma_2:= \sum_{{n_{1}+n_{2}+\cdots+n_{d}=n-1 \atop n_{1},n_{2},\cdots,n_{d}\geq1}
  \atop {0.6s\leq m \leq s-2}}
    {{n} \choose {1, n_{1},n_{2},\cdots,n_{d}}}  2^{\sum_{i=1}^{d}{{n_{i}} \choose 2}}
    \prod_{i=1}^{d-1}(2^{n_{i}}-1)^{n_{i+1}} \nonumber\\
 &\leq
 n_{(d+1)}
  \sum_ {0.6s\leq m \leq s-2} \left(
{{n-d-1} \choose {n_{1}-1,n_{2}-1,\cdots,n_{d}-1}}
 \times
  \exp_2 \left[\sum \binom{n_i}{2} + \sum (n_i n_{i+1}-1)\right] \right)\nonumber\\
 &\leq
 n_{(d+1)} \sum_ {0.6s\leq m \leq s-2} \label{eq:13}
 \Bigg(\bigg( \sum_{m \text{ is fixed}} {{n-d-1} \choose {n_{1}-1,n_{2}-1,\cdots,n_{d}-1}}
 \bigg)\\
 {}&\quad\quad\quad\quad\quad\quad\quad\quad\quad\quad\quad\quad\times
  \exp_2 \left[\max_{m \text{ is fixed}} \left\{ \sum \binom{n_i}{2} + \sum (n_i n_{i+1}-1)
 \right\}\right]\Bigg).  \nonumber
 \end{align}
The total sum of all of the $d$-nomial coefficient of order $s$ is $d^s$,
 the number of $d$-coloring of an $s$-element set.
In the sum (\ref{eq:13}) we add up only those where
 $m=n_{i-1}-1+n_i-1+n_{i+1}-1$ for some
  $2\leq i\leq d-1$.
Choose first an $i$, then $m$ element from the $s$-set, then color those with
  $3$ colors (namely colors $i-1$, $i$ and $i+1$),
 and color the rest by the remaining $d-3$ colors.
We obtain
$$
  \sum_{m \text{ is fixed}} {{n-d-1} \choose {n_{1}-1,n_{2}-1,\cdots,n_{d}-1}}
  \leq (d-2)\binom{s}{m} 3^m (d-3)^{s-m}< (d-2) s^{s-m} 3^s (d-3)^{s-m}.
  $$
Using again (\ref{eq:gf}) and
 then (\ref{eq:m_is_big}) we have
$$
\max_{m \text{ is fixed}} \left\{ \sum \binom{n_i}{2} + \sum (n_i n_{i+1}-1)
 \right\}  \leq \max_{m \text{  is given}}
 f(\mathbf{x}) +\frac{5s}{2} \leq  \frac{1}{2}m^2+ \frac{1}{2}(m-s)^2
 +\frac{5s}{2}.
  $$
So (\ref{eq:13}) gives
$$
  \Sigma_2 \leq n_{(d+1)} \sum_ {0.6s\leq m \leq s-2}
    (d-2) s^{s-m} 3^s (d-3)^{s-m}
   \exp_2 \left[\frac{1}{2}m^2+ \frac{1}{2}(m-s)^2
 +\frac{5s}{2}\right].
  $$
Hence
$$
 \frac{\Sigma_2}{h_1(n,d)}\leq
   \frac{(s+1)(s+2)}{9} 2^s\sum_ {0.6s\leq m \leq s-2}
  \left( s(d-3)2^{-m}\right)^{s-m}.
  $$
One can see that in the given range this sum is dominated by the term
 $m=s-2$, when it is
$O(s^2d^2)2^{-2s+4}$, hence
 $\Sigma_2 = O(s^4d^22^{-s})=o(h_1(n,d))$ follows.
\medskip

\medskip
{\bf Case 3.1.}\quad   $n_{i-1}+n_i+n_{i+1}=s+2$, ($1< i< d$), $n_t=2$ where $t\neq i-2, i+2$,
  and $n_j=1$ for $0\leq j \leq d$, $j\notin \{ i-1, i, i+1, t\}$.

There are $d-2$ ways to choose $i$ then at most $d-3$ possibilities were left to $t$, then
 $n_{(d-3)}$ possibilities to fix $N_j$ $j\neq i-1, i, i+1, t$.
Then one can select $N_t$ in $\binom{s+4}{2}$ ways and distribute the remaining $s+2$
 elements among $N_{i-1}$, $N_i$ and $N_{i+1}$.
Then (\ref{eq:1}) gives

\begin{align}
 \Sigma_{31} &\leq
  n_{(d-3)}(d-2)(d-3)\binom{s+4}{2} \nonumber \\
  {}& \quad\quad\quad \times\sum_{a+b+c=s+2 \atop a,b,c\geq1}
    {{s+2} \choose {a,b,c}}  2^{{a \choose 2}+{b\choose 2}+{c\choose 2}+{2\choose 2}}
   (2^{a}-1)^b (2^b-1)^c (2^c-1)^1(2^2-1)^1\nonumber\\
 &\leq
  12n_{(d-3)}\binom{d-2}{2}\binom{s+4}{2} 2^{\binom{s+2}{2}} \sum_{a+b+c=s+2 \atop a,b,c\geq1}
     {s+2\choose a}{s+2-a\choose c} 2^{-ac+c}.
  \label{eq:c31}
\end{align}
Using standard binomial identities we get
\begin{align}
&\quad \sum_{a+b+c=s+2 \atop a,b,c\geq1}
     {s+2\choose a}{s+2-a\choose c} 2^{-ac+c} \nonumber \\
     &\quad =
     \sum_{a=1, \, 1\leq c<s+1}
    {s+2\choose 1}{{s+1} \choose c}  +
         \sum_{a\geq 2} \binom{s+2}{a}
     \sum_{1\leq c< s+2-a} \binom{s+2-a}{c}(2^{-a+1})^c  \nonumber\\
     &\quad \leq (s+2)2^{s+1}+   \sum_{a\geq 2} \binom{s+2}{a} (1+ 2^{-a+1})^{s+2-a}\nonumber\\
 &\quad \leq  (s+2)2^{s+1} +  \sum_{a\geq 2} \binom{s+2}{a} (3/2)^{s+2-a}
  \leq (s+2)2^{s+1}+ (5/2)^{s+2}.
  \label{eq:c31b}
\end{align}
This is $o(3^s/d)$ so (\ref{eq:c31}) gives $\Sigma_{31}=o(h_1(n,d))$.

The rest of the cases are quite similar.

\medskip
{\bf Case 3.2.}\quad
 $n_{i-1}+n_i+n_{i+1}=s+2$, $n_{i-2}=2$, ($2< i< d$),
 $n_{i+1}\geq 2$,  and $n_j=1$ for $0\leq j \leq d$, $j\notin \{ i-2, i-1, i, i+1\}$.

There are $d-3$ ways to choose $i$, then
 $n_{(d-3)}$ possibilities to fix $N_j$ $j\neq i-2, i-1, i, i+1$.
Then (\ref{eq:1}) gives

\begin{align}
 \Sigma_{32} &\leq
  n_{(d-3)}(d-3) \nonumber \\
    &\quad\quad\quad \times \sum_{a+b+c=s+2 \atop a,b\geq 1,\, c\geq 2}
    {{s+4} \choose {2,a,b,c}}
     2^{{2\choose 2}+{a \choose 2}+{b\choose 2}+{c\choose 2}}
   (2^2-1)^a(2^{a}-1)^b (2^b-1)^c (2^c-1)\nonumber\\
 &\leq
  2 n_{(d-3)}(d-3) \binom{s+4}{2} 2^{\binom{s+2}{2}} \sum_{a+b+c=s+2
 \atop a,b\geq 1, \, c\geq 2}
     {s+2\choose a}{s+2-a\choose c}3^a 2^{-ac+c}.
  \label{eq:c32}
\end{align}
We have
\begin{align}
&\quad \sum_{a+b+c=s+2 \atop a,b\geq 1,\, c \geq 2}
     {s+2\choose a}{s+2-a\choose c}3^a 2^{-ac+c} \nonumber \\
     &\quad \leq
\sum_{a\geq 1,\, 2\leq c < s+2-a}
     {s+2\choose a}{s+2-a\choose c}3^a 2^{-2a+2}
 \leq \sum_{a \geq 1}
     {s+2\choose a} 2^{s+2-a}3^a 2^{-2a+2}
  \nonumber\\
     &\quad = 2^{s+4} \sum_{a\geq 1} \binom{s+2}{a}
   (3/8)^a
 \leq  2^{s+4}(11/8)^{s+2}.
  \label{eq:c32b}
\end{align}
This is $o(3^s)$ so (\ref{eq:c32}) gives $\Sigma_{32}=o(h_1(n,d))$.

\medskip
{\bf Case 3.3.}\quad   $n_{i-1}+n_i+n_{i+1}=s+2$, $n_{i+2}=2$, ($1< i< d-1$),
 $n_{i-1}\geq 3$,  and $n_j=1$ for $0\leq j \leq d$, $j\notin \{ i-1, i, i+1, i+2\}$.

There are $d-3$ ways to choose $i$, then
 $n_{(d-3)}$ possibilities to fix $N_j$ $j\neq i-1, i, i+1, i+2$.
Then (\ref{eq:1}) gives
\begin{align}
 \Sigma_{33} &\leq
  n_{(d-3)}(d-3) \nonumber \\
    &\quad\quad\quad \times \sum_{a+b+c=s+2 \atop a\geq 3,\, b,c\geq1}
    {{s+4} \choose {a,b,c,2}}
     2^{{a \choose 2}+{b\choose 2}+{c\choose 2}+{2\choose 2}}
   (2^{a}-1)^b (2^b-1)^c (2^c-1)^2(2^2-1)\nonumber\\
 &\leq   6n_{(d-3)}(d-3) \binom{s+4}{2} 2^{\binom{s+2}{2}} \sum_{a+b+c=s+2
 \atop a\geq 3, \,b,c\geq1}     {s+2\choose a}{s+2-a\choose c} 2^{-ac+2c}.
  \label{eq:c33}
\end{align}
We have
\begin{align}
&\quad \sum_{a+b+c=s+2 \atop a\geq 3,\, b,c\geq1}
     {s+2\choose a}{s+2-a\choose c} 2^{-ac+2c} \nonumber \\
     &\quad =
         \sum_{a\geq 3} \binom{s+2}{a}
     \left( \sum_{1\leq c< s+2-a} \binom{s+2-a}{c}(2^{-a+2})^c \right) \nonumber\\
     &\quad \leq  \sum_{a\geq 3} \binom{s+2}{a}
 (1+ 2^{-a+2})^{s+2-a}\nonumber\\
 &\quad \leq    \sum_{a\geq 3} \binom{s+2}{a} (3/2)^{s+2-a}
  \leq  (5/2)^{s+2}.
  \label{eq:c33b}
\end{align}
This is $o(3^s)$ so (\ref{eq:c33}) gives $\Sigma_{33}=o(h_1(n,d))$.

\medskip
{\bf Case 4.1.}\quad $n_{i-1}+n_i+n_{i+1}=s+3$,  $n_{i-1}\geq 2$, $n_{i+1}\geq 2$,
   and $n_j=1$ for $0\leq j \leq d$, $j\notin \{ i-1, i, i+1\}$.

There are $d-2$ ways to choose $i$, then
 $n_{(d-2)}$ possibilities to fix $N_j$ $j\neq i-1, i, i+1$.
Then (\ref{eq:1}) gives
\begin{equation}\label{eq:41}
\Sigma_{41} \leq
  n_{(d-2)}(d-2) \times S,
  \end{equation}
where
$$
 S:= \sum_{a+b+c=s+3 \atop a\geq 2,\, b\geq 1, \, c\geq 2}
    {{s+3} \choose {a,b,c}}
     2^{{a \choose 2}+{b\choose 2}+{c\choose 2}}
   (2^{a}-1)^b (2^b-1)^c (2^c-1).
  $$
We separate the case $a=2$ and use obvious upper bounds
\begin{align}
 S &\leq
\sum_{b+c=s+1\atop 2\leq c\leq s}
    {{s+3} \choose 2} {s+1\choose c}
     2^{1+{b\choose 2}+{c\choose 2}}
    3^b 2^{bc+c}
   \nonumber \\
    &\quad\quad\quad
     +
 \sum_{a+b+c=s+3 \atop a\geq 3,\, b\geq 1, \, c\geq 2}
    {{s+3} \choose {a,b,c}}
     2^{{a \choose 2}+{b\choose 2}+{c\choose 2}+ab+bc+c} \nonumber\\
 &=  2{s+3 \choose 2 }2^{\binom{s+1}{2}} \sum_{2\leq c \leq s}
    {s+1\choose c} 3^{s+1-c} 2^c
     \label{eq:41a} \\
    &\quad\quad\quad
     +\quad
  2^{s+3\choose 2}
  \sum_{1\leq b\leq s-2}{s+3 \choose b}
   \left( \sum_{a+c=s+3-b \atop a\geq 3, \, c\geq 2}
    {{s+3-b} \choose {a}}
     2^{-ac+c}\right) .
  \label{eq:41b}
\end{align}
In the row (\ref{eq:41b}), for a given $b$, the terms in the last sum form a unimodal sequence,
 the two terms at the ends are the largest ones.
More precisely, for $ a,c\geq 2$ integers
$$  \frac{{{a+c} \choose {a}}
     2^{-ac+c}}{ {{a+c} \choose {a+1}}
     2^{-(a+1)(c-1)+(c-1)}}=\frac{(a+1)2^{-a}}{c2^{-c}}>1
     \quad \Longleftrightarrow \quad  a\leq c.
  $$
Thus we can upper estimate these terms by the (sum of the) extreme ends,
 when $(a,c)=(3,s-b)$ and when $(a,c)=(s-b+1,2)$.
\begin{align}
 \sum_{a+c=s+3-b \atop a\geq 3, \, c\geq 2}
    {{s+3-b} \choose {a}}
     2^{-ac+c} &\leq (s-1-b)\left(  {{s+3-b} \choose {3}}
     2^{-2s+2b}+  {{s+3-b} \choose {2}}
     2^{-2s+2b}\right) \nonumber \\
     &\leq s^4 4^{-s+b}.
  \nonumber\end{align}
In the row (\ref{eq:41a}) the sum is at most $(3+2)^{s+1}$.
We obtain
\begin{align}
S\leq& (s+3)(s+2)2^{\binom{s+1}{2}}5^{s+1}+ 2^{s+3\choose 2} s^4 4^{-s}
  \sum_{1\leq b\leq s-2}{s+3 \choose b} 4^{b}
  \nonumber\\ \leq& O(s^4)2^{s+1\choose 2}5^s .\nonumber
  \end{align}
This is $o(2^{s+2\choose 2}3^s)$ so (\ref{eq:41}) gives $\Sigma_{41}=o(h_1(n,d))$.

\medskip
{\bf Case 4.2.}\quad $n_{d-1}+n_d=s+2$,  and $n_j=1$ for $0\leq j \leq d-2$.

There are $n_{(d-1)}$ possibilities to fix $N_j$, $j=0,1,\dots, d-2$.
Then (\ref{eq:1}) gives
\begin{align}  
\Sigma_{42} &\leq
  n_{(d-1)}\sum_{a+b=s+2}
    {{s+2} \choose {a}}
     2^{{a \choose 2}+{b\choose 2}}
   (2^{a}-1)^b \nonumber\\
   &\leq n_{(d-1)} \sum   {{s+2} \choose {a}}
     2^{s+2\choose 2}=  n_{(d-1)}2^{s+2\choose 2} 2^{s+2}
      =o(h_1(n,d)). \nonumber
    \end{align}

\medskip
{\bf Case 4.3.}\quad $n_{i-1}+n_i=s+2$, $1< i<d$,
  and $n_j=1$ for $0\leq j \leq d$, $j\notin \{ i-1, i\}$.

There are $d-2$ choices for $i$ and
 $n_{(d-1)}$ possibilities to fix $N_j$, $j=0,1,\dots, d$, $j\neq i-1,i$.
Then (\ref{eq:1}) gives
\begin{align} 
\Sigma_{43} &\leq
  n_{(d-1)}(d-2)\sum_{a+b=s+2}
    {{s+2} \choose {a}}
     2^{{a \choose 2}+{b\choose 2}}
   (2^{a}-1)^b (2^b-1)\nonumber\\
   &\leq n_{(d-1)}(d-2) \sum   {{s+2} \choose {a}}
     2^{s+2\choose 2} 2^b =  n_{(d-1)}(d-2) 2^{s+2\choose 2} 3^{s+2}
      =2h_1(n,d).\nonumber
    \end{align}

\medskip
Adding up the above eight cases, we get that the right hand side
 of (\ref{eq:23}) is at most \hbox{$(2+o(1))h_1(n,d)$},
 completing the proof of the upper bound.
Together with the lower bound~(\ref{eq:2}) we have the asymptotic.

We also obtained that almost all members of ${\mathcal G}(n,d)$
 belong to the group of Case 4.3. One can see that almost all members of the group 4.3. belong to
 ${\mathcal H_1}(n,d)$, thus finishing the proof of Theorem 1.

\section{Upper bound for Theorem \ref{thm2}}\label{upperthm2}

In this section we suppose that $ n-c \log n<d $, where $c$ is a
 sufficiently small constant.
Again we are going to use (\ref{eq:23}).
We put the terms of the right hand side of (\ref{eq:23}) into four groups
 according to $t$, the number of non-singleton classes
$$
  t:=|\{ i: |N_i|>1\}|.
  $$
We have $t\leq n-d-1$.
If $t = n-d-1$, then we have $t$ pairs and $d+1-t$ singletons, i.e.,
 all $n_i\leq 2$.
 \newline
${}$\quad -- Case 1: $t < n-d-1$,  \newline
${}$\quad -- Case 2: $t = n-d-1$ and $\max \{ n_1,n_2, n_{d-2}, n_{d-1}, n_d \} = 2$.\newline
${}$\quad -- Case 3: $t = n-d-1$, $n_d=1$ but there is an
  $i$ with $n_i=n_{i+1}=2$,\newline
${}$\quad -- Case 4: the graphs in ${\mathcal H}_2(n,d)$.  \newline
These exhaust all possibilities.
We will show that the sum in each of the above group is $o(h_2(n,d))$, except
 in the Case 4.
Recall that $2h_2(n,d)= n_{(d+1)}d^s3^s$.

\medskip
{\bf Case 1.}\quad $t < n-d-1:=s$. \newline
Every graph in this class can be obtained by the following five-step
 procedure.\newline
1) Take a path $P:=v_0, v_1, \dots, v_d$, there are $n_{(d+1)}$ ways to do it.
We will have $v_i\in N_i$. \newline
2) Choose $d-t$ indices from $[d]$, the corresponding classes and $v_0$ are the
 singletons, there are $\binom{d}{t}\leq d^t/t!$ ways to do this.\newline
3) Put a second element to the non-singleton classes from the
 $s$ vertices outside the path, there are
$$s_{(t)}=\binom{s}{t} t!=
  \binom{s}{s-t} t!\leq s^{s-t} t!$$
ways to proceed. \newline
4) Distribute the remaining $s-t$ vertices arbitrarily among
 the non-singleton classes, there are $t^{s-t}$ ways of this.
  We now have a partition $(N_0, N_1, \dots, N_d)$ together with a path
  $P$.\newline
5) Finally, call a pair $xy$ {\it open} if either it is contained in some
 $N_i$ or $x\in N_i$, $y\in N_{i+1}$ with $|N_i|>1$ and it is not an edge of $P$.
There are
\begin{equation}\label{eq:51}
  E:=\sum \binom{n_i}{2} + \sum_{n_i>1} n_in_{i+1}-1
  \end{equation}
open pairs.
With given $P$ and a  partition $(N_0, N_1, \dots, N_d)$
 we can select at most $2^E$ subsets of open pairs to create a
 graph from $\mathcal{G} ( x_{0}, N_{1}, \cdots, N_{d})$.

Define $x_i:=n_i-1$ and use  (\ref{eq:gf}) and then (\ref{eq:m_is_small})
 from Lemma~\ref{lm1}.
Note that $m\leq s-(t-3)$, since there are $t$ positive $x_i$'s.
We obtain that the right hand side of (\ref{eq:51}) is at most
$$
 f({\mathbf x})+\frac{5s}{2}\leq \frac{3}{4}(s-t+3)s + \frac{5s}{2}<
  s(s-t)+ 5s.
  $$
So the number of graphs counted in Case 1 is at most
$$
\sum_{1\leq t<s} n_{(d+1)}\times \frac{d^t}{t!}\times
   s^{s-t} t!\times t^{s-t}\times 2^{s(s-t)+5s}
  = 2 h_2(n,d) \left(\frac{32}{3}\right)^s\sum_{s-t\geq 1}
  \left( \frac{st2^s}{d}\right)^{s-t}.
  $$
This is $o(1)$ since the base of the geometric series is $o( (32/3)^{-s})$
 if $s=n-d-1< (\log_2 n)/6$.

\medskip
{\bf Case 2.}\quad $n_j\leq 2$ for all $1\leq j \leq d$, and $\max \{ n_1,n_2, n_{d-2}, n_{d-1}, n_d \} =2$.

We consider the case $n_d=2$ only, the other cases can be handled in the same way.
In this case (\ref{eq:1}) gives at most $2^s9^s$ graphs.
Furthermore there are $\binom{d-1}{s-1}\leq sd^{s-1}/s!$ ways to select the $s$
 indices of the $2$-element blocks.
So the number of partitions with $n_d=2$ is
$$
  \frac{s d^{s-1}}{s!}\times {n \choose 2}{{n-2} \choose 2}
  \cdots { {n -2(s-1)}\choose 2}  (n-2s)!.
  $$
So the number of graphs in this case is at most
$$
  2^s3^{2s}\times\frac{sd^{s-1}}{s!}\frac{n!}{2^s}
  = 2h_2(n,d)\frac{s3^s}{d}
  $$

\medskip
{\bf Case 3.}\quad $n_j\leq 2$, for all $1\leq j \leq d$, $n_d=1$ and
 there is an $i$ with $n_i=n_{i+1}=2$.

Inequality (\ref{eq:1}) gives at most $2^s9^s$ graphs.
Furthermore, there are
$$ \binom{d-1}{s}-\binom{d-s}{s} \leq (s-1)\binom{d-2}{s-1}
  \leq s\binom{d-1}{s-1}\leq \frac{s^2}{d} \frac{d^s}{s!}
$$ ways to select the $s$ indices of the $2$-element blocks from
 $\{ 1, 2, \dots, d-1\}$ such a way that two are next to each other.
So the number of graphs in this case is at most
$$
  2^s3^{2s}\times \frac{s^2}{d} \frac{d^s}{s!}{n \choose 2}{{n-2} \choose 2}
  \cdots { {n -2(s-1)}\choose 2}  (n-2s)!
= 2h_2(n,d)\frac{s^23^s}{d}.
  $$

\medskip
Adding up the above three cases, we get that the number of graphs
 of ${\mathcal G}(n,d)\setminus {\mathcal H}_2(n,d)$
 is at most \hbox{$o(h_2(n,d))$},
 completing the proof of the upper bound in Theorem~\ref{thm2}.

\section{Eccentricity}\label{remarks}

The {\it eccentricity} of a vertex $x$ in the graph $G$
 is the maximum over all vertices of
the length of a shortest path from $x$ to that vertex.
In both Theorems above we in fact proved that an asymptotic for the number of
 $n$-vertex graphs having a vertex of eccentricity $d$.

The error terms in the asymptotics are exponentially small. For
   $3\leq d\leq n- c_1\log n$ we have
\begin{equation}\label{eq:61}
   \frac{   |\mathcal{G}  ( n, {\rm{diam}} =d)|}{h_1(n,d)}=1+O\left(
    d^2s^4(\frac{11}{12})^{s}\right),
  \end{equation}
and for $d> n-c_2\log n$ we have
\begin{equation}\label{eq:62}
   \frac{   |\mathcal{G}  ( n, {\rm{diam}} =d)|}{h_2(n,d)}=1+O\left(
    \frac{s^2 (64/3)^s}{d} \right).
  \end{equation}

\section{ Phase transition}
It would be interesting to investigate the {\bf phase transition}, i.e.,
 the case of $n-d=\Theta(\log n)$.


\begin{thebibliography}{99}


\bibitem{B79}  B. Bollob\'as, {\it Graph Theory}, Springer-Verlag, 1979.


\bibitem{B01}  B. Bollob\'as, {\it Random Graph}, Cambridge Math.~Studies, 2001.


\bibitem{B} D. A. Grable, The diameter of  a random graph with bounded diameter,
 \emph{Random Structure and Algorithm} {\bf 6} (1995), 193--199.


\bibitem{T94} I. Tomescu, On the number of graphs having small diameter,
 \emph{Rev. Roumaine Math. Appl.} {\bf 39} (1994), 171--177.


\bibitem{T96} I. Tomescu, An asymptotic formula for the number
 of graphs having small diameter,
 \emph{Discrete Mathematics} {\bf 156} (1996),  219--228.


\end{thebibliography}
\end{document}